\documentstyle [12pt]{article}

\begin{document}
\begin{center}
\LARGE{\bf The Length of the Shortest Closed Geodesics\\
       on a Positively Curved Manifold}\vspace{5mm}\\
\large Yoe Itokawa \normalsize and \large Ryoichi Kobayashi
\end{center}
\begin{abstract}
We give a metric characterization of the Euclidean sphere in terms of the lower bound
of the sectional curvature and the length of the shortest closed geodesics.
\end{abstract}

\normalsize
\newcommand{\eu}{\mbox{\boldmath{$R$}}}
\newcommand{\z}{\mbox{\boldmath{$Z$}}}
\newcommand{\vnn}{\vspace{5mm}\newline\noindent}
\newcommand{\svnn}{\vspace{1mm}\newline}
\newcommand{\nln}{\newline\noindent}
\newcommand{\n}{\mid}
\newcommand{\wt}{\widetilde}
\newcommand{\ep}{\varepsilon}

\begin{center}\section{Introduction}\end{center}
Let $M$ be a complete connected Riemannian manifold of dimension $d$ and class
$C^{\infty}$.  The study of global structure of closed geodesics on $M$
{\it vis a vis} certain quantitative restrictions on the sectional curvature
$K$ of $M$ has attracted considerable interest.  Henceforth, we assume $k$ to
be a positive constant.  It is well-known that if $K \geq k^2$ on all tangent
2-planes of $M$, then there must exist on $M$ a closed geodesic whose length
is $\leq 2\pi/k$.  The purpose of the present paper is to describe a
rigidity phenomenon observed when this length is extremal on $M$.  More
precisely, we prove
\vnn
{\bf Main Theorem.}\quad{\it If $M$ satisfies $K \geq k^2$ and if the
shortest closed geodesics on $M$ have the length $=2\pi/k$, then $M$ is
isometric to $S^d_k$, the Euclidean sphere of radius $1/k$ in $\eu^{d+1}$.}
\vnn
\indent
Note that we make no assumption about the geodesics' having no
self-intersections.  There exists an example of a 2-dimensional smooth surface
all of whose shortest closed geodesics have self-intersections. 
These examples have some regions where the curvature is negative. 
It is reported that E. Calabi has proved that on a positively curved surface,
at least one of the shortest closed geodesics is always without
self-intersections.\vnn
\indent
We now mention some related rigidity phenomena.  Previously, Sugimoto [Su],
improving on an earlier work of Tsukamoto [Ts], proved
\vnn
{\bf Theorem A.}\quad{\it Suppose that $M$ satisfies $4k^2 \geq K \geq
k^2$.  If $d$ is odd, assume that $M$ is simply connected.  Then, if $M$ has a
closed geodesic of length $2\pi/k$, it is isometric with $S^d_k$.}
\vnn
\indent
Recall that under the curvature assumption of Theorem A, if $M$ is simply
connected, the celebrated Injectivity radius theorem, which is primarily due to
Klingenberg (see [CE (\S\S 5.9,10)], [GKM, \S\S7.5,7] and also [CG],
[KS] and [Sa2]) states that all closed geodesics on $M$ have length
$\geq \pi/k$.
\nln
\indent 
However, we point out that, in general, an assumption on the length of the
shortest closed geodesic is a nontrivially weaker condition than an upper bound
on the sectional curvature.  In fact, it is possible to construct, for any
given $k$ and $\delta$, a Riemannian metric on $S^2$ with $K \geq k^2$ and the
length of the shortest closed geodesic $\delta$-close to $2\pi/k$ but whose
curvature grows arbitrarily large ($S^2$ like surface with highly curved ``equator").  
This construction means that, from the point of view of rigidity theorems in 
Riemannian geometry, imposing an upper bound on the curvature is not natural 
in characterizing a Euclidean sphere among complete Riemannian
manifolds with $K \geq k^2$ having a shortest closed geodesic of
length just $2\pi/k$. 
For more informations on the curvature
bounds and the lengths of closed geodesics, see, for instance, [Sa2].  

In the spacial case of dimension 2, we
have
\vnn
{\bf Theorem B}\,\,({\it Toponogov} [T]).\quad{\it Suppose that $M$
is an abstract surface with Gauss curvature $K \geq k^2$.  If there exists on
$M$ a closed geodesic without self-intersections whose length $=2\pi/k$, then
$M$ is isometric to $S^2_k$.}
\vnn
\indent
However, in higher dimensions, there
are lens spaces of constant sectional curvature $k^2$ so that all geodesics are
closed, the prime ones have no self-intersections, and they are either\nln
\indent
(a) homotopic to 0 and have length  $=2\pi/k$, or\nln
\indent
(b) homotopically nontrivial and can be arbitrarily short.\nln
See [Sa1].  Of course, it follows from our Main Theorem that
\vnn
{\bf Corollary.}\quad{\it If $K \geq k^2$ and the shortest closed
geodesics that are homotopic to 0 in $M$ have the length $2\pi/k$, then the
universal covering of $M$ must be isometric to $S^d_k$.}
\vnn
\indent
Note also that Theorem B is false without the assumption that the closed
geodesics have no self-intersections.  In fact, for any $k$, one can
construct an ellipsoid in $\eu^3$ which possesses a prime closed geodesic 
of length $=2\pi/k$ and whose curvature is $>k^2$.
\nln
\indent
Finally, we mention a previous related result of the first author which gives
another rigidity solution for the nonsimply connected case.
\vnn
{\bf Theorem C}\,\,({\it Itokawa} [I1,2]).\quad{\it If the Ricci
curvature of $M$ is $\geq (d-1)k^2$ and if the shortest closed
geodesics on $M$ have the length $\geq \pi/k$, then either $M$ is
simply connected or else $M$ is isometric with the real projective
space all of whose prime closed geodesics have length $=\pi/k$.}
\vnn
\indent
It is not yet known if our Main Theorem remains true when the assumption on
the sectional curvature is weakened to that on the Ricci curvature.  However, 
we point out that examples were shown in [I1,2] so that for the Ricci curvature
assumption, the shortest closed geodesics may have length arbitrarily close 
to $2\pi/k$ without the manifold's even being homeomorphic to $S^d$.  This
indicates how delicate the Ricci curvarure assumption could be.
\vspace{1mm}
\newline\indent
{\it Acknowledgement.} The authors wish to acknowledge their gratitude 
to D. Gromoll, G. Thorbergsson and W. Ballman for providing them with 
valuable suggestions and informations. 
\vspace{5mm}

\section{Preliminaries}
The purpose of this section is to gather
all the well-known results which will be used in proving the Main Theorem 
as well as to set straight our notational conventions and normalizations.  In
this paper, we agree that by the term {\it curve} we mean an absolutely
continuous mapping $c:\eu \rightarrow M$ whose derivative $\eu \rightarrow
TM$ is an $L_2$ map on each closed interval.  We refer to the restriction of
a curve to any closed interval as an {\it arc}.  If $c$ is a curve and $a<b$
are reals, we write $c_{a,b}$ to denote the arc $c \n_{[a,b]}$.  If $c$
happens to be differentiable, the normal bundle; respectively, the unit
normal bundle of $c$; which are in fact bundles over $\eu$, are denoted
$\perp c$; respectively, $U(\perp c)$.  We shall call a curve $c$ closed if
$c(s+1)=c(s)$ for all $s$.  We denote the set of all closed curves on $M$ by
$\Omega$.\vnn \indent  For fixed $a$, $b$, let $\cal A$$_{a,b}$ denote the
set of all arcs $[a,b] \rightarrow M$.  It is known that $\cal A$$_{a,b}$ has
the structure of a Riemannian Hilbert manifold where the inner product is
given by the natural $L_2$ inner product of variation vector fields along a
curve.  The restriction $c \mapsto c_{0,1}$ embeds $\Omega$ in $\cal
A$$_{0,1}$ as a closed submanifold.  Henceforth, this is the structure we
shall always assume on these spaces.\vnn \indent For $\gamma \in \cal
A$$_{a,b}$, we define the space $\cal V$$'_{\gamma}$ of all square integrable
vector fields $v \in T_{\gamma}\cal A$$_{a,b}$ along $\gamma$ such that
$v(a)=0$, $v(b)=0$, and $v(s) \in \perp_s \gamma$ for all $s$, where $\gamma$
is differentiable.  If $c \in \Omega$, we also define the space $\cal V$$_c$
of all $v \in T_c\Omega$ with $v(s) \in \perp_s c$ almost everywhere.  Then,
$\cal V$$'_{c_{0,1}}$ is canonically embedded in $\cal V$$_c$.\vnn 
\indent We normalize the energy of $\gamma \in \cal A$$_{a,b}$ by\vnn  
$$E(\gamma):=\int_a^b \n \gamma'(s)\n^2ds\,\,. $$
\vnn
Also, we denote by $L(\gamma)$ the length of $\gamma$ in the usual sense.  
Thus, in our convention,
$L(\gamma)^2 \leq E(\gamma)$ with equality if and only if $\gamma$ 
is parametrized proportinal to arclength.  
The term {\it geodesic} is always understood to mean a nonconstant geodesic.  
For $u \in UTM$, the unit tangent bundle, we denote by $c_u$ 
the geodesic $s \mapsto \exp su$.  Recall
that the critical points of $E$ on $\Omega$ are closed geodesics and the
constant curves.\nln \indent
Let $c$ be a geodesic and $a<b \in \eu$.  The Hessian of $E$ at $c_{a,b}$, 
here regarded as a symmetric bilinear form on $T\cal A$$_{a,b}$, is 
denoted $H_a^b$.  We remind that if $v \in \cal V$$'_{c_{a,b}}$ 
and is differentiable outside of finitely many points, 
or if $c \in \Omega$, $a=0$, $b=1$, and $v \in \cal V$$_c$ is
differentiable outside of finite points in $(a,b)$, then $H^b_a(v,v)$ equals
the {\it index integral} (see, for example, [Mi], [Bo] and [BTZ]) 
\begin{eqnarray*}
-2\int_a^b (\langle v''(s),v(s)\rangle + \n v(s)\n^2 \n c'(s)\n^2 K(v(s) \wedge c'(s))ds
-2\sum_s\langle v(s),\Delta_s v'(s)\rangle
\end{eqnarray*}
\noindent
where 
$$\Delta_s v'(s)=v'(s_+)-v'(s_-)$$
denotes the jump in $v'(s)$ at one of its finitely many points of
discontinuity in the open interval $(a,b)$.  We write $\iota'(c_{a,b})$ to
denote the index of $H_a^b\n_{\cal V'}$$_{_{c_{a,b}}}$.  If $c$ is closed, we
put $H:=H_0^1\n_{\cal V}$$_{_c}$ and $\iota(c)$ its index.  We recall the
basic inequality

$$\iota(c) \geq \iota'(c)=\sum_{0<s<1}\nu'(c_{0,s})$$
\noindent
where $\nu'(c_{0,s})$ is the dimension of the space of Jacobi fields in 
$\cal V$$'_{c_{0,s}}$.  In this notation, we state the following well-known
theorem, which is primarily due to Fet [F].\vnn 
{\bf Theorem D.}\quad{\it Assume that
$M$ satisfies $K \geq k^2$.  Then there exists a closed geodesic $c$ on
$M$ such that $L(c) \leq 2\pi/k$ and $\iota(c) \leq d-1$.}\vnn
\indent For each $r \in \eu$, we denote by $\Omega^r$;
respectively, $\Omega^{=r}$ and $\Omega^{<r}$, the subspaces $\{c \in
\Omega\,:\,E(c) \leq r\}$; respectively, $\{c \in \Omega\,:\,
E(c)=r\}$ and $\{c \in \Omega\,;\,E(c) < r\}$.  However, $\Omega^0=\Omega^{=0}$
is identified with $M$ itself and so denoted also by $M$.  It is well-known
that each $\Omega^r$; $r>0$, contains a submanifold $'\Omega_r$ which is
diffeomorphic to an open set in some finite product $M \times
\cdots \times M$ and homotopy equivalent to
$\Omega^{<r}$.  The functional $E$ becomes a proper
function on $'\Omega_r$.  The space $'\Omega_r$ contains all the critical
points in $\Omega^{<r}$ and the Hessian of $E\n_{'\Omega_r}$
retains the same index as $E$ at each critical point.  
For details,
see [Mi (\S 16)] and [Bo].  If $a<r$, we put $'\Omega^a_r:=\,'\Omega_r \cap
\Omega^a$.\nln 
\indent 
We must laler consider a more general
functional $F$ on a finite dimensional Riemannian or separable Hilbert
manifold $X$.  For our purpose, $X$ will be either $\Omega$
or $'\Omega:=\,'\Omega_r$ for some fixed $r$ ($\Omega$ consists of 
absolutely continuous curves with $L_2$ derivative and in particular 
contains contains piecewise differentiable curves, while $'\Omega$ consists 
of broken geodesics).  
Let $c \in X$ be a critical point of $F$.  Then $T_cX$ decomposes
into a direct sum\vspace{1mm}\newline\noindent\begin{center}
$T_cX=$$\cal P$$\oplus$$\cal N$$\oplus$$\cal Z$
\end{center}\noindent
where $\cal P$, $\cal N$ and $\cal Z$ are the spaces on which the Hessian
$H_F$ of $F$ at $c$ is positive definite, negative definite and zero
respectively.  We write $\|\cdot\|$ for the norm in $T_cX$.  Then, we can
state the following important fact due to Gromoll and Meyer [GM1]:\vnn  
{\bf Theorem E}\,\,(Generalized Morse Lemma).\quad{\it In the setting described
above, there exists a neighborhood $U$ of $c$, a coordinate chart
$$\xi_c\,:\,U \longrightarrow T_cX\,\,,$$
\noindent
with respect to which $F$ takes the form
$$F \circ \xi_c^{-1}(v)-F(c) = \| x \|^2-\| y \|^2+f(z)$$
\noindent
where $x$, $y$ and $z$ are the projections of $v \in \xi^{-1}(U)$ on $\cal
P$, $\cal N$ and $\cal Z$ respectively, and the Taylor expansion of $f$ at
$z=0$ starts with a polynomial of degree at least 3 in $z$.  For this
decomposition, $c$ needs not be an isolated critical point of $F$, but if $F$
has other critical points in $U$, they are contained in the zero eigenspace of
$H_F$.}\vnn\indent
We use the notation $U_c^{-0}:=\xi_c^{-1}(\cal N$$\oplus
\cal Z)$ and $U_c^{-}:=\xi_c^{-1}(\cal N)$ and call them the {\it unstable};
respectively, the {\it strong unstable}; submanifolds of $F$ at $c$, even
though we make no assumption on $\dim \cal Z$.
\vnn
\indent 
Suppose that $a \in
\eu$.  We set $\Omega^a_F:=\{c \in \Omega\,\,{\rm or\,\,}'\Omega\,;\,F(c) \leq
a\}$.  Let $I$ be the interval $[-1,1]$.  Suppose $c$ is a critical point of
$F$ with $a:=F(c)$ and $\iota:={\rm index}\,\,H_F\n_{c}=\dim \cal N$.  Let $U$
be a neighborhood of $c$ as defined in Theorem E.  
A differentiable embedding $\sigma:(I^{\iota},\partial I^{\iota}) 
\rightarrow (\Omega,\Omega^a_F-U)$ will be called a {\it weak unstable simplex}
of $F$ at $c$ if $\sigma(0)=c$ and $F\n_{\sigma} \leq a$.  If
$\sigma$ is a weak unstable simplex of $F$ at $c$, then $\sigma \cap
U$ must be contained in the topological cone
$$\{\gamma \in U\,;\,H_F(\xi(\gamma)) \leq 0\}\,\,.$$
\noindent
In particular, if $\sigma(I^{\iota}) \cap \xi^{-1}(\cal N)$ contains an open
neighborhood of $\xi^{-1}(0)$ in $\xi^{-1}(\cal N)$, we call $\sigma$ a {\it
strong unstable simplex}.  We note that
\vnn
{\bf Proposition 1.}\quad{\it Any embedded differentiable simplex
$\sigma:(I^{\iota},\partial I^{\iota}) \rightarrow (X,X-U)$ with $\sigma(0)=c$ and
$H_F(\xi(\sigma(t))) <0$ for all $t \in I^{\iota}-\{0\}$ is homotopic to one
of the strong unstable simplexes modulo $X-U$.}\vnn\indent 
We shall say that a critical point $c$ of $F$ is a {\it nondegenerate}
critical point if $\cal Z$$=\{0\}$.  Note that this agreement is different
from the often-used convention of calling a closed geodesic nondegenerate
if $\cal Z$ is the $S^1$ orbit of the geodesic.  With our convention, a
closed geodesic is never a nondegenerate critical point for $E$.  We put
$a:=F(c)$ and write $X^r:=\{x \in X\,;\,F(x) \leq r\}$.  If $c$ is a
nondegenerate critical point of $F$, then of course $c$ is an isolated
critical point and, for some $\ep>0$, the strong unstable simplexes at $c$
represent a nontrivial class in $\pi_{\iota}(X^{a+\ep},X^{a-\ep})$.\vnn
\indent Let $t>r$.  The set $\Sigma^{\iota}$ of all absolutely continuous
$\sigma:(I^{\iota},\partial I^{\iota}) \rightarrow (X^t,X^r)$ can be given the compact-open
topology.  We will need the following topological\vnn 
{\bf Proposition
2.}\quad{\it Suppose $M$ is compact.  Let $\sigma \in \Sigma^{\iota}$.  Then
there is a neighborhood $\cal O$ of $\sigma$ in $\Sigma^{\iota}$ so that any
$\sigma' \in \cal O$ is homotopic to $\sigma$ modulo $X^r$ up to orientation.}

\section{Proof of the Main Theorem}
It is clear that, in order to prove Theorem 1, we need to consider only one
specific $k$.  So, hereafter we assume that $M$ satisfies $K \geq
k^2$ where $k:=2\pi$.  In the present section, we further assume
that $M$ contains no closed geodesic of length $<1$, or
equivalently that there are no critical points of $E$ in
$\Omega^{<1}-M$.  It now remains for us to prove that then $M$ is
isometric to $S^d_{2\pi}$.\vnn \indent We set\vspace{5mm}
\begin{center}$\cal C$ $:=\{c \in \Omega\,;\,c$ is a closed
geodesic of length 1 and $\iota(c)=d-1\}$ \end{center}\vspace{5mm}
and\vspace{5mm} \begin{center} $\cal C^*$ $:=\{c \in \cal C\,;\,\,$
an unstable simplex of $E$ at $c$ represents\\ a nontrivial element
in $\pi_{d-1}(\Omega,M)\}\,\,.$ \end{center}\vspace{5mm}
Theorem D and the Morse-Schoenberg index comparison assert that $\cal C \not=
\emptyset$.  If we can asume that each $c \in \cal C$ has an isolated 
critical $S^1$-orbit, the technique of Gromoll and Meyer [GM2] fairly readily 
shows that $\cal C^*$ too is non-empty.  In our case, however, it will be precisely 
one of our points that no $c \in \cal C^*$ has an isolated critical orbit. 
Under the stronger hypothesis of $4k^2 \geq K \geq k^2$, Ballman [Ba] showed
that all closed geodesics have nontrivial unstable simplexes.  However, he
makes essential use of the upper bound for $K$ which is not available to us. 
Nonetheless, we shall still prove in the next section,\vnn {\bf Lemma
1.}\quad{\it Under the assumptions of this section, $\cal C^*$ is nonempty
and is a closed set in $\Omega$.}\vnn \indent In this section, we assume
Lemma 1 for the time being, and prove\vnn {\bf Lemma 2.}\quad{\it For each $c
\in \cal C^*$, there is a neighborhood $\cal U$ of $c'(0)$ in $UT_{c(0)}M$
such that whenever $u \in \cal U$ and $\tau$ is any tangent 2-plane
containing $c_u'(s)$ for some $s \in \eu$, then $K(\tau)=k^2$ ($c_u$ being 
the geodesic determined by the initial condition $c_u'(0)=u$).}\vnn \indent
The idea for proving this is to construct for every $c \in \cal C$ a specific
unstable simplex $\sigma$ (the ``Araki simplex") 
and its deformation so that, unless the conclusion of
the lemma is met, $\sigma$ is deformed into $M$, which is a contradiction if $c
\in \cal C^*$.  First, we show\vnn {\bf Lemma 3.}\quad{\it If $c \in \cal C$,
then for any $s \in \eu$ and any $v \in \perp_sc$, $K(c'(s) \wedge v)=k^2$.}\vnn
\indent
{\it Proof.} Assume that, for some $s_1 \in \eu$ and $v_1 \in \perp_{s_1}c$,
$K(c'(s_1) \wedge v_1)>k^2$.  By virtue of the natural $S^1$-action on
$\Omega$, it is no loss of generality to assume that $0 < s_1 < 1/2$.  Now, we
define a real number $\delta$ as follows.  If there is a point in $(0,s_1]$
which is conjugate to 0 along $c$, we choose any $\delta$ so that
$s_1<\frac12-\delta<\frac12$.  If, on the other hand, there is no conjugate
point in $(0,s_1]$, there is a unique Jacobi field $Y$ along $c$ with $Y(0)=0$
and $Y(s_1)=v_1$, and by a consequence of the original Rauch comparison theorem
[CE (\S 1.10, Remark, p.35)], there is an $s_2$, $s_1<s_2<1/2$ so that
$Y(s_2)=0$.  In this case, we choose $\delta$ so that 
$s_2 < \frac12 - \delta < \frac12$.  In either case, we have
$\iota'(c_{0,\frac12-\delta}) \geq 1$.  On the other hand, by the
Morse-Schoenberg index comparison with $S^d_k$, we have
$\iota'(c_{\frac12 - \delta,1}) \geq d-1$, since $L(c_{\frac12 - \delta,1}) \geq
\frac12$.  Therefore, we have
$$\iota(c) \geq \iota'(c) \geq
\iota'(c_{0,\frac12 - \delta}) + \iota'(c_{\frac12 - \delta,1}) \geq 1+d-1=d\,\,,$$
\noindent
which is a contradiction. $\Box$\vnn  
\indent
As a consequence, we have\vnn
{\bf Lemma 4.}\quad{\it Jacobi fields in $\cal V'$$\n_{c_{0,1}}$ are constant
multiple of the fields $\sin(ks)V(s)$; $0 \leq s \leq 1$, where $V$ is any
parallel vector field of elements in $U(\perp c\n_{[0,1)})$ ($\cal V'$$\n_{c_{0,1}}$
being the space of all square integrable normal vector fields along the arc 
$c_{0,1}:=c\n_{[0,1]}$ 
with Dirichlet boundary condition, as is defined in \S2).}\vnn
{\bf Remark.}\quad{\it A priori, the holonomy along the loop $c$ might be
non-trivial.  So, a parallel vector field $V$ of elements in $U(\perp
c\n_{[0,1)})$ might not close up at $s=1$.  Later we will show that the
holonomy along $c$ is trivial.}\vnn 
\indent
Now, let
$V_1,\cdots,V_{d-1}$ be parallel vector fields of orthonormal elements in $U
(\perp c\n_{[0,1)})$.  By compactness argument, there exists an $\eta>0$ so
that each orthogonally trajecting geodesics $t \mapsto \exp tx$ where
$x \in U(\perp c)$ has no point focal to $c$ in $t < \arctan \eta$. 
Define $2(d-1)$ vector fields $X_i(s)$ and $Y_i(s)$ ($0 \leq s \leq 1$) along
$c$ as follows.  These vector fields are not continuous at $s=0$ and
$s=\frac12$. \[ X_i(s)=\left\{ \begin{array}{ll} V_i(s) & {\rm if\,\,}0 \leq s
\leq \frac12\\ 0 & {\rm if\,\,}\frac12<s<1 \end{array}
\right. \]
\noindent
and
\[ Y_i(s)=\left\{ \begin{array}{ll}
0 & {\rm if\,\,}0\leq s \leq \frac12\\
V_i(s) & {\rm if\,\,}\frac12<s<1\,\,.
\end{array}
\right. \]
\noindent
Let $x=(x_1,\cdots,x_{d-1}) \in I \subset \eu^{d-1}$ and $y=(y_1,\cdots,y_{d-1})
\in I \subset \eu^{d-1}$, where $I$ is a small interval in $\eu^{d-1}$.  We define
a $2(d-1)$-dimensional differentiable simplex $\wt\sigma$ 
in $\Omega$ (here we regard $\Omega$ as a Riemannian Hilbert manifold 
consisting of absolutely continuous alosed curves with $L_2$ inner product) 
as follows:
$$\wt\sigma(x,y)(s)=\exp_{c(s)}\arctan\{\eta\sin (2\pi
s)(\sum_{i=1}^{d-1}(x_iX_i(s)+y_iY_i(s)))\}\,\,.$$
\noindent
Here, by  ``arctan" of a vector, we  will mean for a vector $x \in U(\perp c)$ 
to be $(\arctan \| x\|)\frac{x}{\|x\|}$.  W. Ballman pointed out to us that Araki [A]
constructed a simplex in the same way, i.e., varying Jacobi fields
independently outside zero set, when $M$ is a symmetric space.  So, such
a simplex may be called {\it Araki simplex}.  
The Araki simplex $\widetilde\sigma$ consists of curves all passing through
$c(0)$ and $c(\frac12)$.  We deform this simplex in the following way.  If
$x=y$ and if the parallel vector field $\sum_{i=1}^{d-1}x_iV_i$ closes up at
$s=1$, then we make no change on the corresponding loop.  If $x=y$ and if the
paralell vector field $\sum_{i=1}^{d-1}x_iV_i$ does not close up at $s=1$,
then we make a suitable short cut at the non-trivial angle created by the
holonomy.  If $x \not= y$, then we make suitable short cuts at the
non-trivial angle created by the discrepancy $x \not= y$ at $s=\frac12$ and at 
the angle at $s=0$ ($s=1$) if it is non-trivial. 
For instance, we fix a small positive number $\delta$ and make a short cut
between points corresponding to $s=\frac12-\delta$ and $s=\frac12+\delta$. 
After performing this modification and reparametrizing the corresponding
loops by arc length, we get a $2(d-1)$-dimensional simplex $\sigma$ (we call 
this the ``short cut construction"). 
We note that
\svnn

(i) the intersection $\wt\sigma \cap \sigma$ consists
of those closed curves that are generated by $x=y$ where $x(=y)$ satisfies the 
condition that the parallel vector field $\sum_{i=1}^{d-1}x_iV_i$ along $c$ closes 
up at $s=1$, i.e., variations which integrate global Jacobi fields on $c\n_{[0,1)}$, 
\svnn

and\svnn

(ii) the vector
fields $\sin (2\pi s)X_i(s)$ and $\sin(2\pi s)Y_i(s)$ are naturally regarded 
as Jacobi fields along $c\n_{[0,\frac12]}$; respectively, $c\n_{[\frac12,1]}$
which vanish at end points.\svnn

If $x \not= y$, then, after performing the 
above modification, we see that $\sigma(x,y)$ is strictly under the level set
$\Omega^{=1}$ of  $E=1$.  We see this, by applying Rauch type comparison
theorem of Berger (Rauch's second comparison; see [CE (\S 1.10)]) to
variations $$\exp_{c(s)}\arctan\{\sin (2\pi
s)\sum_{i=1}^{d-1}x_iX_i(s)\}\,\,,\qquad s \in [0,\frac12]\,\,,\,\,t \in [0,\eta]$$ 
\noindent of $c\n_{[0,\frac12]}$, and $$\exp_{c(s)}\arctan\{\sin (2\pi
s)\sum_{i=1}^{d-1}x_iY_i(s)\}\,\,,\qquad s \in [\frac12,1]\,\,,\,\,t \in [0,\eta]$$ 
\noindent of $c\n_{[\frac12,1]}$ in $\widetilde\sigma$.  Note that the
corresponding comparison variations in $S_{k}^{2}$ generate great semicircles
fixed at the north and south poles.  
By the same reason, even if $x=y$,
$\sigma(x,y)$ is stricly under the level set $\Omega^{=1}$, if
$\sum_{i=1}^{d-1}x_iV_i$ does not close up at $s=1$.\newline\indent
Summing up, we have\vnn  
{\bf Lemma 5.}\quad{\it There exists a neighborhood $W$ of $c \in
\cal C$ in $\Omega$ so that the $2(d-1)$-dimensional simplex $\sigma \cap W$
is contained in $\Omega^1$.}\vnn \indent We now consider the foliation on
$\eu^{2d-2}$ by affine subspaces orthogonal to the linear subspace defined by
the equations $x=y$.  It follows from Proposition 1 that the $(d-1)$-simplex
$\tau$ which is defined by making suitable short cut modifications to
the collection of piecewise smooth closed curves
$$\tau(x)(s)=\sigma(x,-x)(s)=\exp_{c(s)}\arctan\{\eta\sin(2\pi
s)(\sum_{i=1}^{d-1}(x_iX_i(s)-x_iY_i(s)))\}\,\,,$$ 
\noindent
i.e.,
$\tau(x)=\sigma(x,y)$, where $(x,y)$ varies over the leaf containing the
origin (i.e., $y=-x$), is homotopic to an unstable simplex in the sense of
Proposition 1.  Indeed, we can easily show that, for any vector field $v$ in
$\cal V$$'_{c_{0,1}}$ tangent to $\tau$, $H_0^1(v,v)<0$ for the Hessian
$H_0^1$ at $c$.  This is essentially because of the Taylor expansion $\cos
\theta=1-\frac{\theta^2}{2}+O(\theta^4)$ as $\theta \rightarrow 0$.  Thus,
from Proposition 1 and the consequence of Theorem E, we have
\vnn 
{\bf Lemma 6.}\quad{\it If $c \in \cal
C$$^*$, then there is a neighborhood $U$ of $c$ in $\Omega$ so that, for
$\ep>0$ and a subneighborhood $W \subset U$, $\tau$ constructed above
represents a nontrivial element in $\pi_{d-1}(W,W \cap \Omega^{1-\ep})$.}
\vnn
\indent
The reason why we have chosen the ``short cut construction" is explained 
in the following way. 
If we define a $(d-1)$-dimensional simplex
$\overline \tau$ by 
$$\overline\tau(x):=\exp_{c(s)}\arctan \{\eta\sin(\pi
s)\sum_{i=1}^{d-1}x_iV_i(s)\}\,\,,$$
\noindent then direct calculation of the
Hessian implies that $\overline\tau$ also defines an unstable simplex at $c$
and belongs to the same class as $\tau$ in the relative homotopy group.  This
unstable simplex corresponds to the eigen vector of the index form with
negative eigenvalue.  In this sense, $\overline \tau$ is more natural than
$\tau$.  Now define a $(2d-2)$-simplex $\overline\sigma$, which also contains
the $(d-1)$-simplex (the one defined by $x=y$ in our simplex $\sigma$) 
corresponding to the global Jacobi field on $c$, by 
$$\overline\sigma(x,y)(s):=\exp_{c(s)}\arctan\{\eta\sum_{i=1}^{d-1}(x_i\sin(\pi
s)+y_i\sin(ks))V_i(s)\}\quad (k=2\pi)\,\,.$$
\noindent
Although this construction is natural, it turns out that it is not clear
whether there exists an interval $I$ containing $0$ such that $\overline\sigma(I
\times I)$ is contained in $\Omega^1$.  This is the reason why our
construction of $(2d-2)$-simplex $\sigma$ is based on the short cut argument
of broken geodesics in the model space, although the unstable simplex $\tau$
does not directly integrate the negative eigenspace of the Hessian of the
energy functional $E$ at $c$.\vnn
\indent

We return to our ``short cut construction" and consider the holonomy problem 
mentioned just after Lemma 4. 
One of the following two cases is possible.  Namely, either
\vnn
(A)\,\,For at least one choice of $x_0 \in I$, there is some $\ep>0$ such that
$$\begin{array}{c}
\exp_{c(s)}\arctan\{t\sin(2\pi
s)(\sum_{i=1}^{d-1}x_{0,i}X_i(s)+x_{0,i}Y_i(s))\}\\[5mm]
=\exp_{c(s)}\arctan\{t\sin(2\pi s)(\sum_{i=1}^{d-1}x_{0,i}V_i(s))\}\,\,,
\end{array}$$
\noindent
modified at $s=0$ ($s=1$) if necessary, is contained in $\Omega^{1-2\ep}$ for
all $t \in (\eta/2,\eta]$. In the picture of this situation, we find two 
variation vector fields $V$ and $Y$ along $c$ of  the form 
$$\begin{array}{ccl} Y & = & \sin(2\pi
s)\sum_{i=1}^{d-1}x_iV_i(s)\,\,\,\,\,\,({\rm Jacobi\,\,fields})\,\,,\\[5mm]
V & = & \sin(2\pi s)
\sum_{i=1}^{d-1}(x_iX_i(s)-x_iY_i(s))\,\,\,({\rm tangent\,to}\,\tau)\,\,,
\end{array}$$
outside small neighborhoods of $s=0$ ($s=1$) and $s=\frac12$. 
In $\sigma$, we find a $(d-1)$-dimensional simplex $\tau\cap W$ which 
lies in $\Omega^{<1}$ except at $c$ and moreover we have extra one direction 
represented by $x_0$ which also behaves exactly like a strong unstable 
simplex, or else,
\vnn
(B)\,\,There exists an $\alpha$; $0<\alpha<1$ so that whenever $\n
x_1,\cdots,x_{d-1}\n \leq \alpha$,
$$\exp_{c(s)}\arctan\{\sin(2\pi s)(\sum_{i=1}^{d-1}x_iV_i(s))\}$$
\noindent
lies in $\Omega^{=1}$.  (In particular, each
parallel vector field 
$\sum_{i=1}^{d-1}x_iV_i(s)$ closes up at $s=1$.)
\vnn
\indent 
If Case (A) prevails, $\tau \cap W$ rides on a $d$-dimensional
submanifold of $\sigma \cap W$ which lies in $\Omega^{<1}$ except at $c$ and
hence $\tau \cap W$ can be deformed into $W \cap \Omega^{1-\ep}$.  This
contradicts the conclusion of Lemma 6.  Hence $c \not\in \cal C$$^*$.  If, on
the other hand, we start out with a $c \in \cal C$$^*$, then Case (B) must
really be the case.  In particular, the holonomy along $c \in \cal C$$^*$
must be trivial.  We thus get a $(d-1)$-dimensional local submanifold $S$ of
$\Omega^{=1}$ which is tangent to the 0 eigenspace of the Hessian of $E$
defined on $\cal V$$_{c_{0,1}}$ through $c$.\vnn 
{\bf Lemma 7.}\quad{\it In the present
situation, each parallel vector field
$\sum_{i=1}^{d-1}x_iV_i(s)$ closes up at $s=1$, i.e., the holonomy along $c$
is trivial, and each member $\wt c$ of $S$ is a (smooth) closed geodesic in
$\cal C$.}\vnn\indent
{\it Proof.} We need to prove the second assertion.  If $\wt
c$ is not a critical point of $E$, there exists at least one $x_0 \in I
\subset \eu^{d-1}$ such that the $(d-1)$-dimensional simplex defined by the
$(d-1)$-dimensional affine subspace through $x_0$ orthogonal to the linear
subspace defined by $x=y$ contains no critical point.  Then, by following the
trajectory of $-{\rm grad}\,\,E$, $\tau$ is deformed into $W \cap
\Omega^{1-\ep}$, which contradicts the assumption that we started with $c \in
\cal C$$^*$.  Hence all $\wt c \in S$ are closed geodesics.  If some $\wt c
\not\in \cal C$, then, it follows from Lemma 3 and its proof that $\iota(\wt c)>d-1$. 
So $\tau$ is again deformed into $W
\cap \Omega^{1-\ep}$ via the unstable simplex of $\wt c$.  Hence, either way,
we get a contradiction.  $\Box$\vnn \indent By construction, we also see that
for any $\wt c \in S$, $\wt c(0)=c(0)$.  Translated into $M$, this means that
there is an open tube $B$ (cone-like at $s=0$ and $\frac12$) 
around the set $c(0,\frac12) \cup c(\frac12,1)$
such that for each $q \in B$, a geodesic joining $c(0)$ to $q$ extends to a
closed geodesic in $\cal C$ whose image lies in $B$ except at $s={\rm
half\,\,integers}$.  Applying Lemma 3 to each geodesics proves Lemma 2. 
$\Box$\vnn \indent Even more is true.  By Proposition 2, we get\vnn {\bf
Lemma 8.}\quad{\it Let $c \in \cal C$$^*$ and let $\cal U$$\subset T_{c(0)}M$
be the set in Lemma 2.  Then, there exists an open set $\cal U$$^*$; $c'(0)
\in \cal U$$^* \subset \cal U$, so that, for all $u \in \cal U$$^*$, $c_u \in
\cal C$$^*$.}\vnn \indent
That is to say, the set
\begin{center}
$\cal U$$^*=\{u \in UT_{c(0)}M\,;\,c_u \in \cal C$$^*\}$
\end{center}
is an open set in $UT_{c(0)}M$.  On the other hand, by Lemma 1 and the
continuous dependence of geodesics on their initial values, the set
$\cal U^*$ is also a closed set.  Since $UT_{c(0)}M$ is connected, $\cal U^*$
must in fact be all of $UT_{c(0)}M$.  Together with Lemma 2, we summarize our
result as\vnn
{\bf Lemma 9.}\quad{\it Let $M$ be assumed in this section.  Then, there exists
a point $p \in M$ such that for all $u \in T_pM$, $c_u$ is a closed geodesic of
prime length $1$ and $K(\tau)=k^2$ for all 2-planes $\tau$ tangent to the radial
direction from $p$.}\vnn
\indent
{\it Proof.} Take a $c \in \cal C^*$ and let $p:=c(0)$.  $\Box$\vnn
\indent
Now it is a standard technique to construct an explicit isometry from $M$ onto
$S^d_k$ exactly as in Toponogov's Maximum diameter theorem (see, for instance, 
[CE (\S 6.5)] or [GKM (\S 7.3)]).  Thus, the Main Theorem is proven as soon as
Lemma 1 is established.
\vspace{1mm}

\section{Proof of Lemma 1}
In this section, we continue to assume $K \geq 4\pi^2$.  The following proposition
is essentially contained in some earlier works of M. Berger and is easy to
prove by Morse-Schoenberg index comparison with $S_k^d$ and the tautological
isomorphism $\pi_i(\Omega) \cong \pi_{i+1}(M)$.\vnn
{\bf Proposition
3.}\quad{\it If $M$ contains no closed geodesic of length $\leq 1/2$, then $M$
has the homotopy type of a sphere.  In particular, we have} \[ \pi_i(\Omega,M)
\cong \left \{ \begin{array}{ll} \z & {\rm if\,\,}i=d-1\\
0 & {\rm for\,\,}0 \leq i \leq d-2
\end{array}
\right. \]
\noindent
for the relative homotopy groups $\pi_i(\Omega,M)$ up to $i \leq d-1$.\vnn
\indent
We now return to the assumption that the length of the shortest closed
geodesics on $M$ is 1.  Let $\cal C$ and $\cal C^*$ be as defined in \S 3. 
We wish to prove that a strong unstable simplex at at least one $c \in \cal
C$ represents a nontrivial class of $\pi_{d-1}(\Omega,M)$.  Our technique
will be to approximate $E$ with other functionals that are guaranteed to have
nontrivial unstable simplexes.  Although all our arguments carry through in
all of $\Omega$ in an $S^1$-invariant fashion, essentially because the
functional $E$ satisfies the Condition (C) of Palais and Smale and because an
$S^1$-invariant formulation of Theorem E is available [GM2], we find it a
little easier to work in a finite dimensional space.\nln
\indent
More precisely, choose $r$ sufficiently large, say $r>2$.  Then, all closed
geodesics not in $\Omega^r$ will have index $>2(d-1)$.  Let
$'\Omega:=\,'\Omega_r$.  Then $$\pi_i('\Omega,M) \cong \pi_i(\Omega,M)$$
\noindent
for all $i$; $0 \leq i \leq 2d-3$, and $d-1 \leq 2d-3$ if $d \geq 2$.  Using
Theorem E and a partition of unity on $'\Omega$, we can approximate $E$ with a
sequence $\{E_n\}_{n=1}^{\infty}$ of functionals on $'\Omega$ with the
following properties.\svnn

(i) $\lim_{n \rightarrow \infty} E_n=E$ in the $C^2$
topology.\svnn

(ii) For some $\ep>0$, all critical points of $E_n$ in the
closure of the set $L:=\,{'\Omega_{1+\ep}}-\,{'\Omega_{1-\ep}}$ either belong
to $\Omega_{=1}$ or have index $\geq 2d-2$, and outside $L$, each $E_n$ agrees 
with $E$.\svnn

(iii) Each $E_n$ has only nondegenerate critical points in the
set $L$, all of which have index $\geq d-1$.\svnn

Let $C$ be the set of all closed geodesics in $'\Omega^{=1}$ and let $C_n$ be
the set of all critical points of $E_n$ that lie in $L$.\vnn
{\bf Lemma 10.}\quad{\it For each $n$, there exists in $C_n$, at least one
critical point of $E_n$ that possesses a strong unstable simplex that represents
a nontrivial element in $\pi_{d-1}('\Omega,M)$.}\vnn
\indent
{\it Proof.} From the topology described in Proposition 3, there must exist a
nontrivial element $\rho$ of $\pi_{d-1}('\Omega,M)$.  We first deform $\rho$ so
that the only points of $C_n - (M \cap C_n)$ that lies on the image of $\rho$
are the relative maxima of $E_n \circ \rho$.  In fact, since there are no
critical points of index $<d-1$ except in $M$, at every critical point of $E_n$
lying on $\rho$, say $c$, other than relative maxima, the unstable dimension of
$E_n$ in $'\Omega$ is strictly greater than the unstable dimension of
$E_n \circ \rho$ in the image of $\rho$.  Therefore, in some neighborhood of $c$
in which a chart of the form described in Theorem E is valid, we can deform
$\rho$ in a direction transversal to itself and which decreases $E_n$.  Since
the critical points of $E_n$ are isolated and $\rho$ is contained in a compact
region, by repeating this deformation a finite number of times and by deforming
$\rho$ along the trajectory of $-{\rm grad}\,\,E_n$, we can deform $\rho$ until
it is expressed as a sum of disjoint simplexes, each summand of which is a
simplex in $('\Omega,M)$, hanging from a single critical point of index $=d-1$.
Such critical points must be in $C_n$, and at least one summand must be
nontrivial itself.  $\Box$\vnn
\indent
Of course, it is not necessarily true that a sequence of critical points
$\{c_n\}$ of $C_n$ converges to a closed geodesic.  However, that $\lim_{n
\rightarrow \infty} C_n \subset C$ in the following weaker sense is clear.\vnn
{\bf Lemma 11.}\quad{\it Given any open neighborhood $\cal U$ of $C$ in
$'\Omega$, whenever $n$ is large enough, $C_n \subset \cal U$.}\vnn
\indent
In fact, since the convergence is specified in the $C^2$ topology, we can state
the even stronger\vnn
{\bf Lemma 12.}\quad{\it Let $\{\cal U$$^{-}_c \subset \cal U$$^{-0}_c\}_{c
\in C}$ be a family of pairs of open sets in $'\Omega$ so that, for each
$c \in C$, $\cal U$$^{-}_c$ is a neighborhood of the strong unstable
submanifold $U^{-}_c$ of $E$ at $c$ and $\cal U$$^{-0}_c$ is a neighborhood of
the unstable submanifold $U^{-0}_c$.  Then, for $n$ sufficiently large, for
each $c_n \in C_n$, there exists some $c \in C$, so that $U_{c_n}^{-}$, the
strong unstable manifold of $E_n$ at $c_n$ is contained in $\cal U$$^{-0}_c$. 
Moreover, for such $c_n$ and $c$, a strong unstable simplex $\tau_n$ of $c_n$
contains a subsimplex $\tau'_n$ with $\dim \tau'_n=\dim U_c^-=\iota(c)$ which
is actually contained in $\cal U$$^-_c$.}\vnn
\indent
To see the above, we can take a local coordinate expression around each $c \in
C$ as described in Theorem E and look at the partial derivatives.  By taking
$n$ large, if $c_n \in C_n$ is close to $c \in C$, the corresponding second
derivatives respectively of $E_n$ at $c_n$, $E$ at $c_n$ and $E$ at $c$ can all
be made arbitratily close to each other by the property (i).  But, in $U$, the
strong unstable submanifolds and unstable submanifolds are determined by the
second partial derivatives.\nln
\indent
Now, for each $n$, let $c_n$ be the critical point in Lemma 11 which has a
strong unstable simplex $\tau_n$ that is nontrivial in $\pi_{d-1}(\Omega,M)$. 
For such a $c_n$, $\tau_n \cap U$ must itself be contained in a neighborhood $\cal
U$$^-_c$ of the strong unstable submanifold at some $c \in C$ by index
comparison and the dimensional consideration.  From the construction of
$\tau_n$, this $c$ must be $\in \cal C$.  Let $\tau$ be a strong unstable
simplex at $c$ with $\tau(\partial I^{d-1}) \subset M$.  By repeating the
standard Morse theoretic arguments as in Lemma 6, $\tau$ is seen to
represent a nontrivial element in $\pi_{d-1}(\Omega,M)$.  Hence $c \in \cal
C$$^*$.  Then, that $\cal C$$^*$ is closed follows from Proposition 2.  This
completes the proof of Lemma 1 and thus of Main Theorem.  
$\Box$
\small
Y. I.\nln
Department of Information and Communication Engineering\nln
Fukuoka Institute of Technology\nln
R. K.\nln
Graduate School of Mathematics\nln
Nagoya University
\end{document}